\UseRawInputEncoding
\documentclass[11pt,reqno]{amsart}
\usepackage{fullpage}
\usepackage{mathrsfs,amssymb,graphicx,verbatim,amsmath,amsfonts,mathtools}
\usepackage{paralist}

\usepackage[breaklinks,pdfstartview=FitH]{hyperref}
\usepackage{upgreek}

\hypersetup{
     colorlinks   = true,
     citecolor    = black
}
\hypersetup{linkcolor=black}
\hypersetup{urlcolor=black}

\newcommand{\mnote}[1]{}
\newcommand{\alex}[1]{}
\renewcommand{\le}{\leqslant}
\renewcommand{\ge}{\geqslant}

\renewcommand{\setminus}{\smallsetminus}
\renewcommand{\gamma}{\upgamma}

\renewcommand{\pi}{\uppi}
\newcommand{\Level}{\mathcal{L}}

\newcommand{\T}{\mathsf{T}}

\newcommand{\e}{\varepsilon}

\newtheorem{theorem}{Theorem}
\newtheorem{lemma}[theorem]{Lemma}

\newtheorem{corollary}[theorem]{Corollary}

\newtheorem{conjecture}[theorem]{Conjecture}
\theoremstyle{definition}

\newtheorem{definition}[theorem]{Definition}
\theoremstyle{remark}
\newtheorem{remark}[theorem]{Remark}

\renewcommand{\tau}{\uptau}
\renewcommand{\chi}{\upchi}

\renewcommand{\xi}{\upxi}
\renewcommand{\rho}{\uprho}

\renewcommand{\subset}{\subseteq}

\newcommand{\N}{\mathbb N}

\renewcommand{\theta}{\uptheta}
\renewcommand{\lambda}{\uplambda}

\DeclareMathOperator{\diam}{diam}

\renewcommand{\emptyset}{\varnothing}
\renewcommand{\gamma}{\upgamma}
\renewcommand{\beta}{\upbeta}
\renewcommand{\alpha}{\upalpha}
\renewcommand{\kappa}{\upkappa}
\renewcommand{\psi}{\uppsi}
\renewcommand{\rho}{\uprho}
\renewcommand{\delta}{\updelta}
\renewcommand{\pi}{\uppi}
\renewcommand{\omega}{\upomega}
\renewcommand{\sigma}{\upsigma}

\renewcommand{\eta}{\upeta}

\renewcommand{\kappa}{\upkappa}
\renewcommand{\mu}{\upmu}
\renewcommand{\nu}{\upnu}
\renewcommand{\pi}{\uppi}
\renewcommand{\zeta}{\upzeta}

\newcommand{\bbQ}{\overline{Q^c}}
\newcommand{\bQ}{Q^c}
\newcounter{EnumResume}

\begin{document}

\title{A simple proof of Dvoretzky-type theorem for Hausdorff dimension  in doubling spaces}

\author{Manor Mendel}
\address{(M.~M) Mathematics and Computer Science Department\\ The Open University of Israel\\ 1 University Road\\ P.O. Box 808\\ Raanana 43107\\ Israel}
\email{manorme@openu.ac.il}


\begin{abstract}
The ultrametric skeleton theorem [Mendel, Naor 2013] implies, among other things, the following nonlinear Dvoretzky-type theorem for Hausdorff dimension:
For any $0<\beta<\alpha$, any compact metric space $X$ of Hausdorff dimension $\alpha$  contains
a subset which is biLipschitz equivalent to an ultrametric and has Hausdorff dimension at least $\beta$.
In this note we present a simple proof of the ultrametric skeleton theorem in doubling spaces using Bartal's  Ramsey decompositions [Bartal 2021].
The same general approach is also used to answer a question of Zindulka [Zindulka 2020] about the existence of ``nearly ultrametric" subsets of compact spaces having full Hausdorff dimension.
 \end{abstract}

\maketitle

%

\section{Introduction} \label{sec:intro}

An ultrametric space is a metric space $(U,\rho)$ satisfying the strengthened triangle inequality $\rho(x,y)\le \max\{\rho(x,z),\rho(y,z)\}$ for all $x,y,z\in U$. 
Saying that $(X,d)$ embeds (biLipschitzly) with distortion $D\in [1,\infty)$ into an ultrametric space means that there exists an ultrametric $\rho$ on $X$ satisfying $d(x,y)\le \rho(x,y)\le Dd(x,y)$ for all $x,y\in X$. 
The \emph{ultrametric distortion} of $X$ is the infimum over $D$ for which $X$ embeds in an ultrametric with distortion at most $D$.
Given a metric space $(X,d)$, a point $x\in X$ and a radius $r\in [0,\infty)$, the corresponding closed ball is denoted $B(x,r)=B_d(x,r)=\{y\in X:\ d(y,x)\le r\}$, and the corresponding open ball is denoted $B_d^\circ(x,r)=\{y\in X:\ d(y,x)< r\}$. 
The following theorem was proved in~\cite{MN-hausdorff,MN-skeleton}.

\begin{theorem}
\label{thm:um-skeleton}
For every $\e\in (0,1)$ there exists $C_\e\in (0,\infty)$ with the following property.  Let $(X,d)$ be a compact metric space and let $\mu$ be a Borel probability measure on $X$. Then there exists a closed subset $U\subseteq X$  satisfying:
\begin{itemize}[\quad$\bullet$]
\item The ultrametric distortion of $U$ is $O(1/\e)$.
\item There exists a Borel probability measure $\nu$ supported on $U$ satisfying
\begin{equation}\label{eq:measure growth}
\nu\left(B_d(x,r)\right)\le \left(\mu(B_d(x,C_\e r)\right)^{1-\e},
\end{equation}
for all $x\in X$ and $r\in [0,\infty)$.
\end{itemize}
\end{theorem}

The above theorem, termed  \emph{the ultrametric skeleton theorem} in~\cite{MN-skeleton}, has its roots in 
Dvoretzky-type theorems for finite metric spaces. 
It has applications for algorithms, data-structures, probability, geometric measure theory,
and descriptive set theory. See~\cite{MN-hausdorff,MN-skeleton,KMZ,bartal2021advances}, as well as~\cite[\S8 \& \S9]{Naor-Ribe}, for more details.
One of its corollaries is the following Dvoretzky-type theorem
for Hausdorff dimension. 

\begin{corollary}[\cite{MN-hausdorff}]
\label{cor:dvo-Hausdorff}
For every $\delta\in (0,1)$, every compact metric space $X$ has a closed subset $U\subseteq X$ having Hausdorff dimension  
$\dim_H (U)\ge (1 - \delta)\dim_H(X)$ whose ultrametric distortion
is at most $C/\delta$, for some universal constant $C>0$.
\end{corollary}
As observed in~\cite{MN-hausdorff}, the trade-off between the Hausdorff dimension and the ultrametric distortion of $U$ in Corollary~\ref{cor:dvo-Hausdorff} is asymptotically tight.

\subsection{A simple proof in doubling spaces}
The proof  of Theorem~\ref{thm:um-skeleton} is involved.
The main purpose of this paper is to present a simple  proof
of a special case of Theorem~\ref{thm:um-skeleton}
when $X$ is a \emph{doubling space},
using only a top-down argument on distances.
A metric space $(X,d)$ is called \emph{$\lambda$-doubling} if any
bounded subset $Z\subseteq X$ 
can be covered by at most $\lambda$ subsets of diameter at most 
$\diam_d(Z)/2$. 
It is called \emph{doubling} if it is $\lambda$-doubling for some $\lambda\in\mathbb N$.

\begin{theorem}
\label{thm:um-skeleton-doubling}
Let $(X,d)$ be a $\lambda$-doubling compact metric space, and let $\mu$ be a Borel probability measure on $X$.
Then for every $t\in\mathbb N$ there exists a closed subset $U\subseteq X$  satisfying:
\begin{itemize}[\quad$\bullet$]
\item The ultrametric distortion of $U$ is  at most $8t$.
\item There exists a Borel probability measure $\nu$ supported on $U$ satisfying
\begin{equation}\label{eq:measure growth-doubling}
\nu\left(B_d(x,r)\right)\le \lambda^{2/t} \left(\mu(B_d(x,(16t+1)  r)\right)^{1-1/t},
\end{equation}
for all $x\in X$ and $r\in [0,\infty)$.
\end{itemize}
\end{theorem}

Using Theorem~\ref{thm:um-skeleton-doubling} we obtain a simple proof of Corollary~\ref{cor:dvo-Hausdorff}  in doubling spaces. Both are proved in Section~\ref{sec:simple-proof}. 

Comparing Theorem~\ref{thm:um-skeleton-doubling} to Theorem~\ref{thm:um-skeleton}, 
we see that~\eqref{eq:measure growth-doubling}
is a qualitatively weaker bound than~\eqref{eq:measure growth} 
because the doubling constant appears as a multiplicative factor in the upper bound on $\nu$.
However, the bound in~\cite{MN-hausdorff, MN-skeleton} on the balls' blow-up parameter --- $C_\e$ in \eqref{eq:measure growth} --- 
is only $C_\e=e^{O(1/\e^2)}$, while
the corresponding parameter in~\eqref{eq:measure growth-doubling}
is  $O(1/\e)$. 
Thus, the new proof raises the question of a possible  improvement to Theorem~\ref{thm:um-skeleton}:
\begin{conjecture}
\label{conj:improved-um-skeleton}
Theorem~\ref{thm:um-skeleton} is also true with $C_\e=O(1/\e)$.
\end{conjecture}

In Section~\ref{sec:zindulka} we show that a variant of the proof of Theorem~\ref{thm:um-skeleton-doubling}
can be used to reprove  
a result of Zindulka~\cite{Zindulka}, 
analogous to Corollary~\ref{cor:dvo-Hausdorff},
about
the existence  of  full dimensional subsets of compact metric spaces that are $\beta$-H\"older embedded in ultrametrics, for any $\beta<1$.
The new proof has the advantage of using 
a weaker doubling condition,  solving an open problem from~\cite{Zindulka}.

A more advanced version of the proof of Theorem~\ref{thm:um-skeleton-doubling} has been recently used  in a follow-up 
paper~\cite{mendel2021dvoretzkytype} to obtain a
stronger form of the ultrametric skeleton theorem:
In addition to the upper bound ~\eqref{eq:measure growth-doubling} 
on $\nu$ in terms of $\mu$, 
it also has a lower bound on $\nu$ in terms of $\mu$. 
That theorem is used to prove 
a Dvoretzky-type theorem similar to Corollary~\ref{cor:dvo-Hausdorff} in which both $X$ and $U$ are also \emph{Ahlfors regular}.

\subsection{Bartal's Ramsey decomposition}
For given subsets $P,Q\subset X$ of a metric space $(X,d)$, we denote by $d(P,Q)=\inf\{d(x,y):\; x\in P,\; y\in Q\}$, the lower distance between them. For $x\in X$ we also denote $d(x,Q)=d(\{x\},Q)$.
At the heart of the proof of Theorem~\ref{thm:um-skeleton-doubling} 
is the following elegant lemma 
of Bartal~\cite{bartal2021advances}.
\begin{lemma}[Bartal's Ramsey decomposition lemma~\cite{bartal2021advances}]
\label{lem:bar-ramsey-decomp}
Let $(X,d)$ be a compact metric space and let $\mu$ be a finite Borel measure on $X$. 
For any closed subset $Z\subset X$, $\diam(Z)>0$, 
and integer $t\in\{2,3,\ldots\}$, there exist disjoint closed subsets $P,Q\subset Z$ 
that satisfy
  $d(P,Q)\ge \diam(Z)/(8t)$,
 $\diam(Z\setminus Q)\le \diam(Z)/2$,
and
\begin{equation} \label{eq:bar-ramsey-decomp}
\mu(P)\cdot\mu^*(Z)^{1/t} \ge \mu (Z\setminus Q)\cdot \mu^*(Z\setminus Q)^{1/t}.
\end{equation}
In particular, 
\begin{equation} \label{eq:bar-ramsey-decomp-2}
\frac{\mu(P)}{\mu^*(P)^{1/t}} + \frac{\mu(Q)}{\mu^*(Q)^{1/t}} \ge \frac{\mu(Z)}{\mu^*(Z)^{1/t}}.
\end{equation}
Here we use the convention $0/0=0$, and $\mu^*$ is defined as 
\begin{equation} \label{eq:mu-star}
\mu^*(A)=\sup_{a\in A}\mu(B(a,\diam(A)/4)\cap A).
\end{equation} 
\end{lemma}

Lemma~\ref{lem:bar-ramsey-decomp} is stated slightly differently in~\cite{bartal2021advances}  and only for finite spaces. 
For completeness, we prove it in Appendix~\ref{sec:proof-bartal}.
Informally, Theorem~\ref{thm:um-skeleton-doubling} is deduced from
Lemma~\ref{lem:bar-ramsey-decomp} as follows.
Lemma~\ref{lem:bar-ramsey-decomp} is applied  iteratively, obtaining a hierarchy of  decompositions that defines an ultrametric on a subset $U\subseteq X$,
with the following property:
Consider a cover of an ultrametric ball $B\subset U$ by ultrametric balls $B_1,B_2,\ldots,B_m \subset U$. 
Those ultrametric balls correspond
to subsets $Q$ and $P_1,P_2\ldots,P_m$ in the hierarchy of decompositions (respectively). A recursive application
of~\eqref{eq:bar-ramsey-decomp-2} gives
\[
\sum_i\frac{\mu(P_i)}{\mu^*(P_i)^{1/t}} \ge \frac{\mu(Q)}{\mu^*(Q)^{1/t}}.
\]
Associate with each ultrametric ball $B$ the value
$\xi(B)=\mu(Q)/ \mu^*(Q)^{1/t}$ 
(recall that $Q\subseteq X$ is the subset corresponding to $B$).
The above inequality means that $\xi$ is sub-additive on the balls of the ultrametric, and hence a ``sub-measure".
In doubling spaces one can relate 
$\mu\asymp\mu^*$, and therefore
$\xi(B)\asymp \mu^{1-1/t}(Q)$. 
Using an argument from~\cite{MN-skeleton}, 
the sub-measure $\xi$ can be transformed into an actual measure $\nu$ on $T$ which is dominated by $\mu^{1-1/t}$
on balls of the metric $d$.

Bartal's Ramsey decomposition lemma and its variants found many applications in metric Ramsey theory. 
See~\cite{bartal2021advances} for a survey of these applications. 
This paper adds another one, 
though it is still unclear whether Bartal's lemma, or some variant thereof, can be used to prove Theorem~\ref{thm:um-skeleton} 
(or Conjecture~\ref{conj:improved-um-skeleton}).

\section{Proof of Theorem~\ref{thm:um-skeleton-doubling}}
\label{sec:simple-proof}

The approximate ultrametric subset extracted in Theorem~\ref{thm:um-skeleton-doubling} has a natural hierarchical structure that we represent using a binary tree as follows.

\begin{definition}[Binary trees]
\label{def:btree}
Let $\{0,1\}^{<\omega}$ be the set of finite binary sequences, and
let $|\cdot|:\{0,1\}^{<\omega} \to \mathbb N$ be the length of sequences.
Let $\preceq$ be the partial order of extension over $\{0,1\}^{<\omega}$, i.e., for $u,v\in \{0,1\}^{<\omega}$
$u\preceq v$ if $|u|\le |v|$, and $u_i=v_i$ for every $i\in\{1,\ldots,|u|\}$.
For $u\in\{0,1\}^k$, we denote
$u^\smallfrown0,u^\smallfrown1\in\{0,1\}^{k+1}$ the unique elements in $\{0,1\}^{k+1}$ for which
$u\preceq u^\smallfrown0$, and 
$u\preceq u^\smallfrown1$.

A (rooted, full, possibly infinite) binary tree $\T\subset \{0,1\}^{<\omega}$ is a non-empty subset satisfying:
\begin{compactitem}
\item For every $u\preceq v\in \T$, we have $u\in\T$; 
\item $u^\smallfrown 0\in \T$ if and only if $u^\smallfrown 1\in \T$.
\end{compactitem}
The elements of $\T$ are called \emph{vertices}.
The vertex $\emptyset\in\{0,1\}^{<\omega}$  must be in every binary tree and is called the \emph{root}.
 For $a\in\{0,1\}$, 
if $u\in\T$ is a vertex of $\T$ and $u^\smallfrown a\in \T$ is also a vertex of $\T$ then $u^\smallfrown a$ 
is called a child of $u$ (in $\T$).
A vertex $u\in \T$ without children in $\T$ is called a leaf (in $\T$).

A subset $b\subseteq \T$ of a binary tree is called a branch in $\T$ if $b$ is  a binary tree in which every vertex has at most one child in $b$, and if this child is a leaf in $b$ then it is also a leaf in $\T$.
The set of branches is called the boundary of $\T$ and is denoted
$\partial_\T(\emptyset)$. 
More generally for a vertex $u\in \T$,
the boundary of $u$ is defined as
$\partial_\T(u)=\{b\in\partial_\T(\emptyset):\; b\upharpoonright |u|=u\}$.

The least common ancestor of a pair of distinct branches $b,c\in\partial_\T(\emptyset)$, denoted $u=b\wedge c\in\T$, is the unique vertex such that $u_i=b_i=c_i$, for $i\in\{0,\ldots,|u|\}$,
and $b_{|u|+1}\ne c_{|u|+1}$. The least common ancestor of two vertices is defined similarly.
\end{definition}

\begin{lemma} \label{lem:compact-um}
Let
$\T$ be a binary  tree. 
Let $\Delta:\T\to [0,\infty)$ be labels on the vertices of $\T$
that are monotone, i.e.,   $\Delta(v)\le \Delta(u)$ for every 
$u\preceq v\in\T$, $\Delta(v)=0$ if and only if $v$ is a leaf, 
and for every branch $b\in\partial_\T(\emptyset)$,
$\inf_{v\in b}\Delta(v)=0$.
Define a distance $\rho$ on $x,y\in \partial_\T(\emptyset)$ as follows:
$\rho(x,y)=\Delta(x\wedge y)$.
Then $(\partial_\T(\emptyset),\rho)$ is a compact ultrametric.
Furthermore,
$\mathcal{O}_\T=\{\partial_\T(v): v=\emptyset \;\vee\; \forall u\prec v, \;(\Delta(v)<\Delta(u)) \}$ 
is the set of open balls in $(\partial_\T,\rho)$, and the set of closed balls in $(\partial_\T,\rho)$ with positive radii.
\end{lemma}

Lemma~\ref{lem:compact-um} is straightforward to prove and folklore. Versions of it can be found, e.g., in~\cite[{\S}2]{MN-skeleton}.
For completeness, we provide a proof in Appendix~\ref{sec:proof-bartal}.

\begin{lemma}
Let $(X,d,\mu)$ be compact metric measure space whose doubling constant is $\lambda\in\N$.
Then, for any Borel set $A\subset X$, 
\begin{equation}\label{eq:mu-mu*}
 \mu^*(A) \le \mu(A)\le \lambda^2 \mu^*(A),
\end{equation}
where $\mu^*$ is defined in~\eqref{eq:mu-star}.
\end{lemma}
\begin{proof}
Let $\Delta=\diam_d(A)$.
By the doubling condition there exists $\lambda^2$ subsets
$A_1,\ldots A_{\lambda^2}\subseteq A$, 
each of diameter at most $\Delta/4$ that covers $A$. 
Let $x_i\in A_i$ be chosen arbitrarily.
Since $B(x_i,\Delta/4)\cap A\supseteq A_i$. The collection of sets
$\{B(x_i,\Delta/4)\cap A\}_{i=1,\ldots \lambda^2}$ covers $A$, and therefore at least one of them satisfies
$\mu(B(x_i,\Delta/4)\cap A)\ge \mu(A) /\lambda^2$. Since the supremum
in~\eqref{eq:mu-star} ranges over all those sets, we conclude that
$\mu^*(A)\ge \mu(A)/\lambda^2$.

In the other direction,
\(
\mu^*(A) =\sup_{a\in A}\mu(B(a,\diam(A)/4)\cap A)
\le \sup_{a\in A} \mu (A)=\mu(A). 
\)
\end{proof}

\begin{lemma}
\label{lem:um-skeleton-doubling}
Fix a compact $\lambda$-\emph{doubling} metric space $(X,d)$, 
a Borel probability measure $\mu$ on $X$, and  $t\in\{2,3,4,\ldots \}$.
Then there exists a binary
tree $\T$, $\Delta:\T\to [0,\infty)$ with the following properties.
Associated with every $u\in \T$ is a ``cluster" $C_u\subseteq X$ satisfying:
\begin{enumerate}[{\quad}(A)]
\item \label{it:root} $C_\emptyset=X$.
\item \label{it:compact} $C_u$ is a closed subset for every $u\in \T$.
\item \label{it:laminar} $C_v\subseteq C_{u}$ for every $u\preceq v\in\T$.
\item \label{it:Delta} $\Delta(u)= \diam_d(C_u)$ 
for every $u\in \T$.
\item \label{it:separation} If $u, v\in\T$, $u\wedge v \notin\{u,v\}$,  then $d(C_u,C_v)\ge \Delta({u\wedge v})/(8t)$.
\item \label{it:singleton} For every branch $b\in \partial_\T(\emptyset)$, $\bigcap_{v\in b} C_v$ is a singleton.  
\setcounter{EnumResume}{\value{enumi}}
\end{enumerate}
By item~\eqref{it:singleton}, we can define a mapping
$\imath:\partial_\T(\emptyset)\to X$, by the set-equation
$\{\imath(b)\}=\bigcap_{v\in b} C_v$.  Then
\begin{enumerate}[{\quad}(A)]
\setcounter{enumi}{\value{EnumResume}}
\item \label{it:partialT=T} 
$\imath(\partial_\T(u))\subseteq C_u$ for every $u\in \T$.
\item \label{it:injective} The mapping $\imath$ is injective.
\setcounter{EnumResume}{\value{enumi}}
\end{enumerate}
Defining 
$U=\imath(\partial_\T(\emptyset))\subset X$, and a distance  $\rho(x,y)=\Delta(\imath^{-1}(x)\wedge \imath^{-1}(y))$ for  $x,y\in U$, then:
\begin{enumerate}[{\quad}(A)]
\setcounter{enumi}{\value{EnumResume}}
\item \label{it:rho}
$(U,\rho)$ is a compact ultrametric and
$d(x,y) \le \rho(x,y)\le 8t \cdot d(x,y)$.
\setcounter{EnumResume}{\value{enumi}}
\end{enumerate}
Lastly, there exists a  function $\xi:\T\to[0,\infty)$
satisfying:
\begin{enumerate}[{\quad}(A)]
\setcounter{enumi}{\value{EnumResume}}
\item \label{it:xi subadditive} $\xi$ is sub-additive:
for every non-leaf vertex $u\in T$,
$\xi(u)\le \xi(u^\smallfrown 0)+\xi(u^\smallfrown 1)$.

\item  \label{it:xi-mu} 
For every $u\in T$, 
\begin{equation}
\label{eq:xi-mu}
\mu^{1-1/t}(C_u)\le \xi(u) \le \lambda^{2/t} \mu^{1-1/t}(C_u).
\end{equation}
\end{enumerate}
\end{lemma}

\begin{proof}
The tree $\T$ and the clusters associated with it are defined recursively. 
The cluster associated with  $\emptyset\in\T$ is $C_\emptyset=X$, which satisfies item~\eqref{it:root}.
Assume that a vertex $u\in\T$ and the associated cluster $C_u$ were defined.
If $C_u$ is a singleton then $u$ will remain without children, i.e., a leaf. 

Assume next that $C_u$ is not a singleton.
By application of Lemma~\ref{lem:bar-ramsey-decomp} on $Z=C_u$ we have $P,Q\subset C_u$ satisfying 
Lemma~\ref{lem:bar-ramsey-decomp}.
Define $C_{u^\smallfrown0}=P$ and $C_{u^\smallfrown 1}=Q$. 
This, in particular,
satisfies items~\eqref{it:compact}, and~\eqref{it:laminar}.
Define $\Delta(u)=\diam_d(C_u)$ which satisfies item~\eqref{it:Delta}.

Item~\eqref{it:separation}.
Suppose $u,v\in\T$ such that $u\wedge v \notin\{u,v\}$
(i.e., $u\not\preceq v$ and $v\not\preceq u$).
Let $w=u\wedge v$.   Assume, without loss of generality that $w^\smallfrown 0\preceq u$, and $w^\smallfrown 1\preceq v$.
 By Lemma~\ref{lem:bar-ramsey-decomp},
 \[
 d(C_u,C_v)\ge d(C_{w^\smallfrown 0}, C_{w^\smallfrown 1}) \ge \Delta(u\wedge v)/(8t),\]
 which proves item~\eqref{it:separation}.

Item~\eqref{it:singleton}.  
First assume that the branch $b$ is finite. In this case $b\in \T$ is a leaf, and $\cap_{v\in b} C_v=C_b$ which is a singleton by the construction above. 
Next assume that $b$ is infinite.
Denote $v_i= b\upharpoonright i$. 
By the above construction, all vertices in $b$ have two children, so
let $u_i$ be the ``other child" of $v_i$, i.e., the child of $v_i$ such that $u_i\ne v_{i+1}$.
$(\Delta({v_i}))_i$ is a positive, non-increasing sequence. 
Assume towards a contradiction that 
$\lim_{i\to\infty} \Delta(v_i)>0$.
That is, there exists $\e>0$ such that $\Delta({v_i})\ge \e$ for every $i\in\mathbb N$. 
Observe that for $i<j$, $u_i \wedge u_j=v_i$.
Therefore,
for every $i\ne j$,
\[
d(C_{u_i},C_{u_j})\ge \Delta({v_{\min\{i,j\}}})/(8t)\ge \e/(8t),
\]
which contradicts the compactness of $X$.
Hence, $\lim_{i\to\infty}\diam_d(C_{v_i})=\lim_{i\to\infty} \Delta({v_i})=0$. By Cantor's intersection theorem,
 $\bigcap_{v\in b} C_v$ is a singleton.

Item~\eqref{it:partialT=T} is an immediate corollary of
the definition of $\imath$ and Item~\eqref{it:singleton}.  
By Item~\eqref{it:separation}, different branches are  mapped to different points of $X$ by $\imath$, 
so Item~\eqref{it:injective} holds.

Item~\eqref{it:rho}. 
It follows from Lemma~\ref{lem:compact-um} that  $(U,\rho)$ is a compact ultrametric.
Fix $x,y\in U$, and denote $w=\imath^{-1}(x)\wedge \imath^{-1}(y)$.
Since $x,y\in C_{w}$,
$d(x,y)\le \diam_d(C_w)=\rho(x,y)$.
On the other hand, since $x\in C_{w^\smallfrown 0}$ while 
$y\in  C_{w^\smallfrown 1}$ (or vice-versa),
\[ 
d(x,y) \ge d(C_{w^\smallfrown 0}, C_{w^\smallfrown 1}) \ge \Delta(x\wedge y)/(8t)=
\rho(x,y)/(8t).
\]
This proves Item~\eqref{it:rho}.

Define $\xi:\T\to[0,\infty)$  as
\[
\xi(u)=\frac{\mu(C_u)}{\mu^*(C_u)^{1/t}}.
\]
By inequality~\eqref{eq:bar-ramsey-decomp-2} of Lemma~\ref{lem:bar-ramsey-decomp}, 
$\xi$ is sub-additive om $\T$,
i.e., $\xi(u)\le \xi(u^\smallfrown 0)+\xi(u^\smallfrown 1)$ for every non-leaf vertex $u\in\T$.
This proves Item~\eqref{it:xi subadditive}.
Lastly, Inequality~\eqref{eq:xi-mu} (Item~\eqref{it:xi-mu}) 
follows directly from \eqref{eq:mu-mu*}.
\end{proof}

\begin{proof}[Proof of Theorem~\ref{thm:um-skeleton-doubling}]
The case $t=1$ is trivial, so we assume from here on that $t\in\{2,3,4,\ldots\}$.
Apply Lemma~\ref{lem:um-skeleton-doubling}
on the metric space $(X,d)$, and the Borel probability measure $\mu$.
We will define the measure $\nu$ supported on $U$
and check that it satisfies the properties asserted
by the theorem.

First, define $\nu: \{\imath(\partial_\T(v))\}_{v\in \T}\to [0,1]$ recursively on $\T$ as follows.
Define $\nu(U)=1$, and note that
$\nu(\imath(\partial_\T(\emptyset)))=1=\mu(X)^{1-1/t}\le \xi(\emptyset)$. 
Assume now that $\nu(\imath(\partial_T(u)))$ was defined,
and define
\[
\nu(\imath(\partial_\T({u^\smallfrown 0}))) = \frac{\xi(u^\smallfrown 0)}{\xi(u^\smallfrown 0)+
\xi(u^\smallfrown 1)} \cdot\nu(\partial_\T(u)), \qquad 
\nu(\partial_\T({u^\smallfrown 1})) = \frac{\xi(u^\smallfrown 1)}{\xi(u^\smallfrown 0)+\xi(u^\smallfrown 1)} \cdot\nu(\partial_\T(u)).
\]
From the inductive hypothesis and the sub-additivity of $\xi$,
\[
\nu(\imath(\partial_\T({u^\smallfrown 0}))) \le 
{\xi(u^\smallfrown 0)} \cdot\frac{\xi(u)}{\xi(u^\smallfrown 0)+\xi(u^\smallfrown 1)}
\le \xi(u^\smallfrown 0),
\]
and in similar fashion $\nu(\imath(\partial_\T({u^\smallfrown 1 })))\le\xi(u^\smallfrown 1)$.
Furthermore, $\nu$ is additive on $\{\imath(\partial_\T(v))\}_{v\in \T}$.

Recall that a (set theoretic) semi-ring in $X$
is a collection of subsets $\mathcal S\subseteq 2^X$ satisfying
\begin{inparaenum}[(i)]
\item $\emptyset\in \mathcal S$;
\item if $A,B\in\mathcal S$ then $A\cap B\in\mathcal S$;
\item if $A,B\in\mathcal S$ then there exist $n\ge 0$ and
$A_1,\ldots, A_n\in\mathcal S$ pairwise disjoint such that 
$A\setminus B=\bigcup_{i=1}^n A_i$.
\end{inparaenum}
By Lemma~\ref{lem:compact-um},
 $\{\emptyset\}\cup \{\partial_\T(v)\}_{v\in\T}$ is 
a semi-ring consisting of all the open balls in $(U,\rho)$,
and
$\nu$ is a pre-measure on that semi-ring. 
By Carath\'eodory extension theorem (see~\cite[Theorem~1.53]{Klenke}), $\nu$ can be extended to a measure on the $\sigma$-algebra generated by $\{\partial_\T(v)\}_{v\in \T}$, 
which is the $\sigma$-algebra of Borel sets of $(U,\rho)$.
Since the metrics $d|_U$ and $\rho$ are topologically equivalent, and $U$ is closed and hence Borel set of $(X,d)$, $\nu$ can be extended to a Borel measure on $X$ by simply define
$\nu(A)=\nu(A\cap U)$ on every Borel $A\subseteq X$.

We are left to prove~\eqref{eq:measure growth-doubling}.
Fix $x\in X$ and $r\ge 0$.
If $B_d(x,r)\cap U=\emptyset$, then
$\nu(B_d(x,r))=0$ and there is nothing to prove.
Otherwise, let $y\in B_d(x,r)\cap U$, so
\[
B_d(x,r)\cap U \subseteq B_d(y,2r) \cap U \subseteq B_\rho(y,16tr).
\]
Since $B_\rho(y,16tr)$ is a closed ball in $(U,\rho)$, by Lemma~\ref{lem:compact-um}  there exists
some $v\in\T$ such that $B_\rho(y,16tr)=\imath(\partial_\T(v))$.
In particular $\diam_d(C_v)=\Delta(v)\le 16tr$.
Observe that 
\[
B_\rho(y,16tr)=\imath (\partial_\T(v))\subseteq C_v \subseteq B_d(y,16tr) \subseteq B_d(x,(16t+1)r).
\]
Hence,
\begin{equation*}
\nu (B_d(x,r))\le \nu(\imath(\partial_\T(v))) \le \xi(v) 
\stackrel{\eqref{eq:xi-mu}}{\le}\lambda^{2/t} \mu(C_v)^{1-1/t}
\le\lambda^{2/t}
  {\mu(B_d(x,(16t+1)r))^{1-1/t}}. \qedhere
\end{equation*}
\end{proof}

\begin{proof}[Proof of Corollary~\ref{cor:dvo-Hausdorff} in doubling spaces]
We follow the argument in~\cite{MN-skeleton}.
Fix $\e\in(0,1/2)$, and let
$t=\lceil 1/\e \rceil$.
Assume that $(X,d)$ is a complete doubling metric space of Hausdorff dimension greater than $\alpha\in (0,\infty)$, and doubling constant $\lambda\in \mathbb N$. 
Then there exists~\cite{How95,Mattila} an $\alpha$-Frostman measure on $(X,d)$.
That is, a positive Borel  measure $\mu$ satisfying $\mu(B_d(x,r))\le Kr^\alpha$ for every $x\in X$ and $r\in (0,\infty)$, where $K>0$ is a constant that may depend  on $X$ and $\alpha$ but not on $x$ and $r$. 

By taking $R>0$ sufficiently large, we may assume that
$\mu(B_d(x,R))>0$. So we restrict ourselves to $X=B_d(x,R)$, and by normalizing $\mu$, we may assume that
$\mu(X)=1$. Since a closed bounded subset of complete doubling space is compact, we may assume that $X$ is compact.

An application of Theorem~\ref{thm:um-skeleton-doubling} to $(X,d,\mu)$ yields a compact subset $U\subseteq X$ that embeds into an ultrametric space with distortion $8t\le 12/\e$, 
and a Borel probability measure $\nu$ supported on $U$ satisfying 
\[
\nu(B_d(x,r))\le \lambda^{2/t} \mu(B_d(x,(16t+1) r))^{1-1/t}\le \lambda^{2/t}K^{1-1/t}(16t+1)^{(1-1/t)\alpha}r^{(1-1/t)\alpha},
\]
 for all $x\in X$ and $r\in (0,\infty)$. 
 Hence, $\nu$ is a $((1-1/t)\alpha)$-Frostman measure on $(U,d)$, implying (see~\cite{Mattila}) that $(U,d)$ has Hausdorff dimension at least $(1-1/t)\alpha \ge (1-\e)\alpha$. 
\end{proof}

\section{Remarks}
\label{sec:remarks}

The following remarks were communicated to the author by an anonymous referee.

In~\cite{MZ}, an analog of  Corollary~\ref{cor:dvo-Hausdorff} for the packing dimension was proved. 
Their proof is an analog of the corresponding proof for Hausdorff dimension in~\cite{MN-hausdorff}, see~\cite[\S2,\S3]{MZ}.
In particular a key step in their proof is the use
of Theorem~\ref{thm:um-skeleton}. 
Therefore, Theorem~\ref{thm:um-skeleton-doubling} also gives a simpler proof of their main result for doubling spaces.

It is possible to obtain a result qualitatively similar 
to Corollary~\ref{cor:dvo-Hausdorff} in doubling spaces using Assouad's embedding~\cite{Assouad}. 
However, this approach, which we outline next,
is quantitatively far from being asymptotically tight (unlike the approach used in this paper).
Let $(X,d)$ be a $\lambda$-doubling compact metric space having Hausdorff dimension $\dim_H(X,d)=\alpha$. Then $\dim_H(X,\sqrt{d})=2\alpha$.
By Naor-Neiman's quantitative version of  Assouad's embedding, see~\cite{Naor-Neiman},
there is a biLipschitz embedding $f:(X,\sqrt{d})\to \mathbb R^n$ with distortion $O(\log \lambda)$ and $n=O(\log \lambda)$. Since biLipschitz embeddings do not change the Hausdorff dimension, $\dim_H(f(X))=2\alpha$.
Consider the following Cantor set $C^n_\e\subseteq \mathbb R^n$:
\[
C^n_\e=\Bigl \{\tfrac{1+\e}{2}\sum_{i=0}^\infty \bigl(\tfrac{1-\e}{2}\bigr)^i \delta_i:\; \delta_i\in\{0,1\}^n\Bigr \}.
\]
The ultrametric distortion of $C_\e^n$ is at most $\sqrt{n}/\e$, and its Hausdorff dimension is $\dim_H(C_\e^n)=n+\log_2(1-\e)$.
By~\cite[Theorem~13.11]{Mattila}, there exists 
an isometry $\phi:\mathbb R^n \to \mathbb R^n$,
such that
\[
\dim_H(\phi(C_\e^n)\cap f(X))\ge \dim_H(C_\e^n)+\dim_H(f(X))-n=2\alpha+ \log_2(1-\e).
\]
Let $S=f^{-1}(\phi(C_\e^n)\cap f(X))\subseteq X$. Since Hausdorff dimension is an invariant of biLipschitz embedding, 
$(S,\sqrt{d})$ is a closed subset of $(X,\sqrt{d})$ 
whose Hausdorff dimension is at least
$2\alpha+(\log_2(1-\e))$ and its ultrametric distortion is at most order of
$\log \lambda \cdot \sqrt{n}/\e$.
Observe that square of the distances of an ultrametric is still an ultrametric, and therefore
$(S,{d})$ is a closed subset of $(X,{d})$ 
whose Hausdorff dimension is at least
$\alpha+(\log_2(1-\e))/2$ and its ultrametric distortion is at most order of
$\log^2 \lambda \cdot {n}/\e^2 \lesssim(\log^3 \lambda)/\e^2$.
By setting $\delta =-  \log_2(1-\e)/(2\alpha)$, we have in terms of Corollary~\ref{cor:dvo-Hausdorff}'s notation, a subspace $(S,d)\subset (X,d)$ whose Hausdorff dimension is at least $(1-\delta)\alpha$, 
and its ultrametric distortion is $O(\log^3\lambda/(\alpha^2 \delta^2))$. Since $\alpha\lesssim \log \lambda$, this bound is no better than $O((\log \lambda) /\delta^2)$.

\section{Nearly ultrametric subsets of full dimension}
\label{sec:zindulka}

In~\cite{Zindulka}, Zindulka proved a qualitative variant of Corollary~\ref{cor:dvo-Hausdorff} with a different 
trade-off for spaces with a mild doubling condition: 
The subset $U\subseteq X$ extracted is full dimensional, i.e.,
$\dim_H(U)=\dim_H(X)$. 
But $U$ is not quite biLipschitz equivalent to an ultrametric, rather it is only ``nearly biLipschitz equivalent".
Rigorously:

\begin{definition}[Nearly ultrametric space~{\cite[Def.~2.4, Prop.~2.5]{Zindulka}}]
A bounded metric space $(U,d)$ is called \emph{nearly ultrametric} if there exists an ultrametric $\rho$ on $U$ such that%
\begin{itemize}[$\bullet$ ]
\item The identity map $(U,\rho)\mapsto (U,d)$ is $1$-Lipschitz. I.e., for every $a,b\in U$, $d(a,b)\le \rho(a,b)$.
\item The identity map $(U,d)\mapsto (U,\rho)$ is 
\emph{nearly Lipschitz}, i.e., for every $\beta\in(0,1)$  there exists $C>0$ such that  
\( \rho(a,b)\le C\cdot d(a,b)^\beta\),
for every $a,b\in U$.
\end{itemize}
\end{definition}

\begin {definition}[Modulus of the doubling condition]
Let $(X,d)$ be  a metric space. 
For $\delta>0$,  $\lambda_X(\delta)$ is defined to be the minimal number $\lambda\in \mathbb N\cup\{\infty\}$ such that any subset of diameter at most $\delta$ can be covered by $\lambda$ subsets of diameter at most $\delta/16$.%
\footnote{The constant $16$ is somewhat arbitrary, and any other constant larger than $1$ can be chosen without altering the meaning of conditions~\eqref{eq:non-exploding} and ~\eqref{eq:non-exploding-2}.  
It was chosen  so that $\lambda_X(\delta)$ would be ``compatible" with $\tilde \lambda_Z$ in~\eqref{eq:AP-decomp} and thus spare the need to 
introduce more notation.}
\end{definition}

See~\cite{Zindulka} for more context, motivation and applications of those definitions.

\begin{theorem}[{\cite[Theorem~2.6]{Zindulka}}]
\label{thm:zindulka}
Let $(X,d)$ be a compact metric space%
\footnote{The original theorem applies more generally to analytic spaces. 
However, the gist of the argument is in the compact case, 
and the analytic case can be reduced to the compact case, see~\cite{Mercourakis,Zindulka}.}
whose modulus of the doubling condition satisfies%
\footnote{In \cite[Def.~2.1]{Zindulka}  this condition is called
\emph{non-exploding}. 
There, $\lambda_X(\delta)$ is replaced with the function $Q(\delta)$ which is defined as ``[the minimal value such that] every closed ball in $X$ of radius $\delta > 0$ is covered by at most $Q(\delta)$ many closed balls of radius $\delta/2$."
These two variants of non-exploding are equivalent, as evident from the (straight-forward) estimates
$Q(\delta/2)\le \lambda_X(\delta) \le  \prod_{i=0}^5 Q(\delta/2^i)$.}
\begin{equation}
\label{eq:non-exploding}
\lim_{\delta\to 0^+}  \frac{\log \lambda_X(\delta)}{\log \delta}=0,
\end{equation}
and let $\mu$ be a finite Borel measure on $X$.
Then for every $\e>0$ 
there exists a closed subset $U\subset X$ for which
$\mu(X\setminus U)<\e$, and $U$ is nearly ultrametric.
\end{theorem}

The following corollary of Theorem~\ref{thm:zindulka}  should be compared with Corollary~\ref{cor:dvo-Hausdorff}.
\begin{corollary} \label{cor:nearly-UM-full-dimension}
Let $(X,d)$ be a compact metric space satisfying~\eqref{eq:non-exploding}
and assume that
$\mathcal H^\alpha(X)>0$, where $\mathcal H^\alpha$ is the $\alpha$-dimensional Hausdorff measure on $X$
(In particular, $\dim_H(X)\ge \alpha$.) 
Then $X$ contains a nearly ultrametric closed subset $U\subseteq X$
having $\dim_H(U)=\alpha$.
\end{corollary}

Zindulka's original proof of Theorem~\ref{thm:zindulka} shares the same high-level approach described in Section~\ref{sec:remarks}:
Apply an Assouad-type embedding (this time, into an infinite dimensional torus) and intersects it with Cantor-like subset.
It was suggested by an anonymous referee that the general approach in this paper may give a different proof of Theorem~\ref{thm:zindulka}. 
It turns out that this is indeed true.
In fact, we obtain the conclusion of Theorem~\ref{thm:zindulka} under a milder doubling condition, requiring only that $\lambda_X(\delta)=e^{\delta^{-o(1)}}$
instead of $\lambda_X(\delta)=\delta^{-o(1)}$ as in~\eqref{eq:non-exploding}. Formally:
\begin{theorem}
\label{thm:zindulka-2}
Let $(X,d)$ be a compact
metric space whose modulus of the doubling condition satisfies
\begin{equation}
\label{eq:non-exploding-2}
\lim_{\delta\to 0^+}  \frac{\log \log (e\lambda_X(\delta))}{\log \delta}=0,
\end{equation}
and let $\mu$ be a finite Borel measure on $X$.
Then for every $\e>0$ 
there exists a closed subset $U\subset X$ for which
$\mu(X\setminus U)<\e$, and $U$ is nearly ultrametric.
\end{theorem}

In particular, Corollary~\ref{cor:nearly-UM-full-dimension}  also holds assuming~\eqref{eq:non-exploding-2} instead of~\eqref{eq:non-exploding}. 
Theorem~\ref{thm:zindulka-2} also answers in affirmative Question~5.9 from~\cite{Zindulka}.
We next outline the proof of Theorem~\ref{thm:zindulka-2}.
We begin with an analog to Bartal's Ramsey decomposition lemma.

\begin{lemma} \label{lem:AP-decomp}
Let $(X,d)$ be a compact metric space and let $\mu$ be a finite Borel measure on $X$. 
For any closed subset $Z\subset X$, $\diam(Z)>0$, 
and integer $t\in\{2,3,\ldots\}$, there exist disjoint closed subsets $P,Q\subset Z$ 
that satisfy
  $d(P,Q)\ge \diam(Z)/(8t)$,
 $\diam(Z\setminus Q)\le \diam(Z)/2$,
and
\begin{equation} \label{eq:AP-decomp}
\mu(P)\cdot\tilde \lambda_Z^{1/t} \ge \mu (Z\setminus Q).
\end{equation}
Here 
$\tilde \lambda_Z$ is the minimal number of subsets of diameter at most $\diam(Z)/8$ needed to cover $Z$.
\end{lemma}

Lemma~\ref{lem:AP-decomp} is a variant of Awerbuch and Peleg's \emph{sparse partitions}~\cite{AP}. 
Its proof is similar to (and simpler than) 
the proof of Lemma~\ref{lem:bar-ramsey-decomp}.
We prove both of them simultaneously in Appendix~\ref{sec:proof-bartal}. 
Next is an analog of Lemma~\ref{lem:um-skeleton-doubling}.

\begin{figure}[ht]
\includegraphics{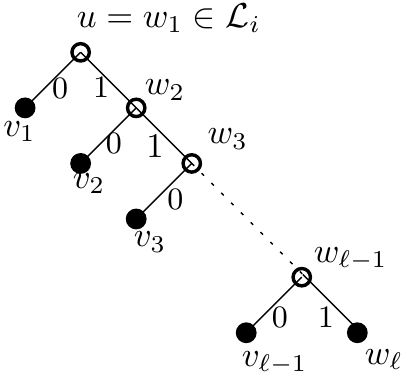}
\caption{The vertices $v_1,\ldots, v_{\ell-1},w_\ell$ constitutes $\Level_{i+1}\cap\T_u$.}
\label{fig:L_i+1}
\end{figure}

\begin{lemma}
\label{lem:nearly-um-skeleton}
Fix a compact metric space $(X,d)$, a Borel probability measure $\mu$ on $X$, and 
$t:\mathbb N \to \{2,3,4,\ldots\}$ non-decreasing.
Then there exists a binary
tree $\T$, $\Delta:\T\to [0,\infty)$ with the following properties.
Associated with every $u\in \T$ is a ``cluster" $C_u\subseteq X$ satisfying:
\begin{enumerate}[{\quad}(A)]
\item \label{it2:root} $C_\emptyset=X$.
\item \label{it2:compact} $C_u$ is closed for every ${u\in\T}$.
\item \label{it2:laminar} $C_v\subseteq C_{u}$ for every $u\preceq v\in\T$.
\item \label{it2:Delta} $\Delta(u)= \diam_d(C_u)$ for every $u\in\T$.
\item \label{it2:separation} If $u, v\in\T$, $u\wedge v \notin\{u,v\}$,  then $d(C_u,C_v)\ge \Delta({u\wedge v})/\bigl(8\cdot t\bigl (\lfloor - \log_2 \Delta(u\wedge v)\rfloor \bigr )\bigr)$.
\item \label{it2:singleton} For every branch $b\in \partial_\T(\emptyset)$, $\bigcap_{v\in b} C_v$ is a singleton.  
\setcounter{EnumResume}{\value{enumi}}
\end{enumerate}
By item~\eqref{it2:singleton}, we can define a mapping
$\imath:\partial_\T(\emptyset)\to X$, by the set-equation
$\{\imath(b)\}=\bigcap_{v\in b} C_v$.  Then
\begin{enumerate}[{\quad}(A)]
\setcounter{enumi}{\value{EnumResume}}
\item \label{it2:partialT=T} 
$\imath(\partial_\T(u))\subseteq C_u$ for every $u\in \T$.
\item \label{it2:injective} The mapping $\imath$ is injective.
\setcounter{EnumResume}{\value{enumi}}
\end{enumerate}
Denote
$U=\imath(\partial_\T(\emptyset))\subset X$, and a distance  $\rho(x,y)=\Delta(\imath^{-1}(x)\wedge \imath^{-1}(y))$ for  $x,y\in U$. Then:
\begin{enumerate}[{\quad}(A)]
\setcounter{enumi}{\value{EnumResume}}
\item \label{it2:rho}
$(U,\rho)$ is a compact ultrametric and
\begin{equation} \label{eq:nearly-um-bound}
d(x,y) \le \rho(x,y)\le 8 \cdot t \bigl(\lfloor - \log_2 d(x,y)\rfloor \bigr) \cdot d(x,y).
\end{equation}
\setcounter{EnumResume}{\value{enumi}}
\end{enumerate}
Denote $\Level_i=\{u\in\T:\; \Delta(u)\le 2^{-i} \land (\forall v\prec u,\; \Delta(v)> 2^{-i})\}$.
Then:
\begin{enumerate}[{\quad}(A)]
\setcounter{enumi}{\value{EnumResume}}
\item \label{it2:mu-decrease}
For every $i\in\mathbb N$,
\begin{equation}
\label{eq:mu(C_u)-bound}
\mu\Bigl(\bigcup _{u\in\Level_i} C_u\Bigr)
\le \lambda_X(2^{-i})^{1/t({i})}
\mu\Bigl(\bigcup _{u\in\Level_{i+1}} C_u\Bigr).
\end{equation}
\setcounter{EnumResume}{\value{enumi}}
\end{enumerate}
\end{lemma}
\begin{proof}[Outline of the proof]

The construction of the tree $\T$ and the labels
$\Delta:\T\to [0,\infty)$ is similar to the construction in the proof of 
Lemma~\ref{lem:um-skeleton-doubling} but instead of applying
Lemma~\ref{lem:bar-ramsey-decomp} to $C_u$ where $u\in\T$,
we apply Lemma~\ref{lem:AP-decomp} to $C_u$ with the parameter $t(\lfloor - \log_2 \Delta(u)\rfloor )$.
Items~\eqref{it2:root}, \eqref{it2:compact}, \eqref{it2:laminar},
\eqref{it:Delta}, \eqref{it2:separation},
\eqref{it2:singleton}, \eqref{it:partialT=T},
 \eqref{it2:injective}, and~\eqref{it2:rho}
are all proved exactly the same as in the proof of Lemma~\ref{lem:um-skeleton-doubling}.

We are left to prove Item~\eqref{it2:mu-decrease}.
For $u\in \T$, denote by $\T_u$ the subtree of $\T$ containing $u$ and all its descendants. To prove~\eqref{eq:mu(C_u)-bound} it is sufficient to prove
that for every $u\in \Level_i$,
\begin{equation} \label{eq:L_i<=L_i+1}
\mu(C_u)\le \lambda_X (2^{-i})^{1/t({i})} \sum_{v\in \Level_{i+1}\cap \T_u} \mu(C_v).
\end{equation}
We fix $u\in\Level_i$.
If $u\in \Level_{i+1}$, then~\eqref{eq:L_i<=L_i+1} holds trivially, so we may assume that $u\notin\Level_{i+1}$.
Since $u\notin\Level_{i+1}$, $u$ is not a leaf 
(observe that for a leaf vertex $w$, $\Delta(w)=0$, and therefore if $w\in  \Level_{i}$ then $w\in\Level_j$ for any $j\ge i$).
We claim that the vertices of $\Level_{i+1}\cap \T_u$ have the structure depicted in Fig.~\ref{fig:L_i+1}:
Let $w_1=u$, and $w_{k+1}=w_k{}^\smallfrown 1$,
let $\ell$ be the smallest $k$ such that $\diam(C_{w_k})\le 2^{-i-1}$. Then $w_\ell \in\Level_{i+1}\cap \T_u$,
and for $k<\ell$,
$v_k=w_k{}^\smallfrown 0\in \Level_{i+1}\cap \T_u$, since, by Lemma~\ref{lem:AP-decomp},
\[
\Delta(v_k)= \Delta(w_k{}^\smallfrown 0)=\diam_d(C_{w_k{}^\smallfrown 0})
\le \diam_d(C_{w_k})/2 \le \Delta(u)/2\le 2^{-i-1}.
\]
Since $\Level_{i+1}$ cannot contain two vertices with ancestor/descendant relation, we deduce that
$\Level_{i+1}\cap \T_u=\{v_1,\ldots v_{\ell-1},w_\ell\}$.
Applying Inequality~\eqref{eq:AP-decomp} inductively, and observing that $\tilde \lambda_{C_{w_k}}\le  \lambda_X(2^{-i})$, 
for $k\in \{1,\ldots,\ell-1\}$,
we conclude that
\begin{multline*}
\mu(C_u)\le  \lambda_X(2^{-i}))^{1/t(i)}\mu(C_{v_1}) + \mu(C_{w_2})
\le \ldots \le \sum_{k=1}^{\ell-1} \lambda_X(2^{-i})^{1/t(i)}
\mu(C_{v_k})+ \mu(C_{w_\ell})\\
\le 
\lambda_X(2^{-i})^{1/t({i})} \Bigl(\sum_{k=1}^{\ell-1} \mu(C_{v_k})+ \mu(C_{w_\ell})\Bigr).
\qedhere
\end{multline*}
\end{proof}

\begin{proof}[Proof of Theorem~\ref{thm:zindulka-2}]
The assertion is trivial when $\mu(X)=0$, so we assume that $\mu(X)>0$.
By rescaling the measure we may actually assume without loss of generality that $\mu$ is a Borel probability measure, i.e., $\mu(X)=1$.
Assume also for now that $\diam_d(X)\le 1/2$.
We will get rid of this condition at the end of the proof.
Fix $\e\in(0,1)$, and define 
\begin{align}\label{eq:def:eta-t}
\eta(i)&=\frac{\log_2 \log (e\lambda_X(2^{-i}))}{i} & t(i)= &\Bigl\lceil \e^{-1}2^{i(\eta(i)+ ({2\log_2 e})/{\sqrt{i}}\,)} \Bigr\rceil.
\end{align}
By~\eqref{eq:non-exploding-2},
\( 
\lim_{i\to\infty}  \eta(i)=0.
\)

Apply Lemma~\ref{lem:nearly-um-skeleton} to $(X,d)$ and  $\mu$ with $t(\cdot)$.
Inequality~\eqref{eq:nearly-um-bound} implies that $(U,d)$ is nearly ultrametric: $(U,\rho)$ is an ultrametric and dominates the distances in $(U,d)$.
In the other direction, fix $\beta<1$.
Since $\lim_{i\to \infty} \eta(i)+({2\log_2 e})/{\sqrt{i}}=0$, there exists
some $\xi\in \mathbb N$ such that for any $i>\xi$,
\(
{{1-\eta(i)-{(2\log_2 e)}/{\sqrt{i}}}}\ge \beta.
\)
By analyzing two  cases: $d(x,y)\ge 2^{-\xi}$, for which we use $\rho(x,y)\le 1/2$, and  $d(x,y)< 2^{-\xi}$, for which we set $i=\lfloor - \log_2
 d(x,y)\rfloor$ and use~\eqref{eq:nearly-um-bound},
we deduce  that
\begin{multline*}
\rho(x,y)\le \max\bigl\{2^\xi,
 8 \cdot  t(i)
 \bigr\} \cdot d(x,y) \\
 \le  \max\Bigl \{2^\xi, 16 \e^{-1}2^{\eta(i)+(2\log_2 e)/\sqrt{i}}
 d(x,y)^{-\eta(i)-(2\log_2 e)/\sqrt{i}}\Bigr\} \cdot d(x,y) \\
\le  \max\{2^{\xi},32\e^{-1}\} d(x,y)^\beta.
\end{multline*}

We next bound $\mu(U)$ from below.
Observe that 
\(U= \bigcap_{i\in\mathbb N} \bigcup _{u\in\Level_{i}} C_u\), and 
\((\bigcup _{u\in\Level_{i}} C_u)_i\) 
is a non-increasing sequence of subsets. 
Therefore,
\(
\lim_{i\to\infty} \mu \bigl(\bigcup _{u\in\Level_{i}} C_u \bigr)=\mu(U).
\) 
By~\eqref{eq:mu(C_u)-bound},
\[
\prod_{j=1}^\infty \lambda_X(2^{-j})^{-1/t(j)} 
\le \prod_{j=1}^\infty \frac{\mu\Bigl(\bigcup_{u\in\Level_{j+1}} C_u\Bigr)}{\mu\Bigl(\bigcup_{u\in\Level_{j}} C_u\Bigr)} =\frac{\mu(U)}{\mu\Bigl(\bigcup_{u\in\Level_{1}} C_u\Bigr)}=\mu(U).
\]

Estimating the infinite product above on the left-hand side:
\begin{equation*}
\prod_{j=1}^\infty \lambda_X(2^{-j})^{1/t(j)}
\stackrel{\eqref{eq:def:eta-t}}\le \prod_{j=1}^\infty \exp\Bigl(2^{j\eta(j)}
\e 2^{-j\eta(j)} e^{-2\sqrt{j}} \Bigr)
= \prod_{j=0}^\infty \exp\Bigl(
\e e^{-2\sqrt{j}}\Bigr)
=   \exp\biggl(\e \sum_{j=1}^\infty
e^{-2\sqrt{j}}\biggr)
< e^{\e/2},
\end{equation*}
where the last inequality follows from the estimate
\[
\sum_{j=1}^\infty e^{-2\sqrt{j}}\le e^{-2}+ \int_1^\infty e^{-2\sqrt{x}}dx= e^{-2}+ \int_2^\infty 0.5 t e^{-t}dt = e^{-2}+1.5 e^{-2}< 0.5.
\]
In summary,
\(
\mu(U)>e^{-\e/2}\ge 1-\e/2.
\)

When $\diam_d(X)>1/2$, we rescale the metric as follows.
Define $\alpha=1/(2\diam(X))$ and the metric 
$\tilde d(x,y)=\alpha\cdot d(x,y)$, so $\diam_{\tilde d}(X)=1/2$.
By the above, there exists $U\subseteq X$ such that
$\mu(U)>1-\e/2$, and an ultrametric $\tilde \rho$ on $U$ 
for which  for any $\beta<1$ there exists $C>0$ such that
\[ \tilde d(x,y)\le \tilde \rho(x,y)\le C\cdot \tilde d(x,y)^\beta.
\]
Define $\rho=\alpha^{-1}\tilde \rho$, which is also an ultrametric on $U$,
and
\[  d(x,y)\le  \rho(x,y)\le C \alpha^{-1}\cdot \tilde d(x,y)^\beta= C\alpha^{\beta-1}\cdot d(x,y)^\beta. 
\qedhere
\]
\end{proof}


\appendix

\section{ } 
\label{sec:proof-bartal}

Bartal's proof of Lemma~\ref{lem:bar-ramsey-decomp} is written for finite spaces, but applies to compact spaces with only minor adaptations. 
We include the proof here for completeness.
We also include a proof of  Lemma~\ref{lem:AP-decomp}.
Lemma~\ref{lem:AP-decomp} can have a slightly simpler proof (with slightly better parameters), but for brevity we deduce it from the proof of Lemma~\ref{lem:bar-ramsey-decomp}.

\begin{proof}[Proof Lemma~\ref{lem:bar-ramsey-decomp} and Lemma~\ref{lem:AP-decomp}]
Assume $\infty>\mu(Z)>0$ (otherwise, not much to prove).
Denote $\Delta=\diam(Z)$.
With the convention $0/0=0$, let $x_P\in Z$ be a point that maximizes
\[
\frac{\mu (B(x,\Delta/8)\cap Z)} {\mu(B^o(x,\Delta/4)\cap Z)}.
\]
With this choice, $\mu (B(x_P,\Delta/8)\cap Z)>0$.
For $i\in\{0,1,\ldots,t-1\}$, let $H_i=B(x,(1+i/t)\Delta/8)\cap Z$, and also define
$H_t=B^o(x_P,\Delta/4)\cap Z$.
Clearly there exists $i\in\{1,\ldots,t\}$ for which
\begin{equation}\label{eq:light-ring}
\mu(H_i)\le \mu(H_{i-1}) \cdot \biggl(\frac{\mu(H_t)}{\mu(H_0)}\biggr)^{1/t}.
\end{equation}
We then set $P=H_{i-1}$,  $\bQ=B^o(x_P,(1+i/t)\Delta/8)$,  $Q=Z\setminus \bQ$. 
Clearly, 
\[\diam(\bQ)\le \diam(B^o(2\Delta/8))\le \Delta/2.\]
Observe that 
$P$ and $Q$ are compact, 
and
$H_i\supseteq \bQ\supseteq P$.
By the triangle inequality, for every $a\in P$, and $b\in Q$, $d(a,b)\ge d(b,x_P)-d(a,x_P)\ge \Delta/(8t)$.
The measure satisfies
\begin{equation} \label{eq:light-ring-2}
\mu(P) \stackrel{\eqref{eq:light-ring}}{\ge} 
\mu(\bQ) \cdot \biggl(\frac{\mu(B(x_P,\Delta/8)\cap Z)}{\mu(B^o(x_P,\Delta/4)\cap Z)}\biggr)^{1/t}.
\end{equation}

Let $u\in \bbQ= B(x_P,(1+i/t)\Delta/8 )$ the point that maximizes
$\mu(B(u,\diam(\bbQ)/4)\cap \bbQ)$.
Since $\diam(\bbQ)\le \Delta/2$,
we have
$B(u,\diam(\bbQ)/4)\cap \bbQ\subseteq B(u, \Delta/8 )\cap Z$, and hence
$\mu^*(\bbQ)\le \mu(B(u,\Delta/8)\cap Z)$.
Also,
$\mu^*(Z)\ge \mu(B^o(u,\Delta/4)\cap Z)$.
From the definition of $x$
we therefore have
\[
\frac{\mu(B(x_P,\Delta/8)\cap Z)}{\mu(B^o(x_P,\Delta/4)\cap Z)}
\ge 
\frac{\mu(B(u,\Delta/8)\cap Z)}{\mu(B^o(u,\Delta/4)\cap Z)}
\ge 
\frac{\mu^*(\bbQ)}{\mu^*(Z)} 
\ge \frac{\mu^*(\bQ)}{\mu^*(Z)} . 
\]
Plugging the last inequality into~\eqref{eq:light-ring-2},
we obtain~\eqref{eq:bar-ramsey-decomp}. 
Inequality~\eqref{eq:bar-ramsey-decomp-2} follows 
from~\eqref{eq:bar-ramsey-decomp}:
\begin{multline*}
\frac{\mu(P)}{\mu^*(P)^{1/t}} + \frac{\mu(Q)}{\mu^*(Q)^{1/t}}
\stackrel{\eqref{eq:bar-ramsey-decomp}}{\ge}
\frac{\mu(\bQ)}{\mu^*(P)^{1/t}} \cdot \frac{\mu^*(\bQ)^{1/t}}{\mu^*(Z)^{1/t}}
+  \frac{\mu(Q)}{\mu^*(Q)^{1/t}}
 \ge \frac{\mu(\bQ)}{\mu^*(Z)^{1/t}} +\frac{\mu(Q)}{\mu^*(Z)^{1/t}}
= \frac{\mu(Z)}{\mu^*(Z)^{1/t}}. 
\end{multline*}

We next show that the subsets $P$ and $Q$ also satisfy~\eqref{eq:AP-decomp}, thus proving Lemma~\ref{lem:AP-decomp}.
Indeed, by the definition of $\tilde\lambda_Z$ and $x_P$,
\[
\frac{\mu(B(x_P,\Delta/8)\cap Z)}{\mu(B^o(x_P,\Delta/4)\cap Z)} =
 \max_{x\in Z}\frac{\mu(B(x,\Delta/8)\cap Z)}{\mu(B^o(x,\Delta/4)\cap Z)} \ge 
 \max_{x\in Z}\frac{\mu(B(x,\Delta/8)\cap Z)}{\mu(Z)} \ge 
 \tilde \lambda_Z^{-1}.
\]
Plugging it into~\eqref{eq:light-ring-2}, we obtain
\(
\mu(P) \cdot \tilde \lambda_Z^{1/t}\ge \mu(Q^c). 
\)
\end{proof}

\begin{remark}
In the proof above we assumed that 
\[ f(x)=\frac{\mu (B(x,\Delta/8)\cap Z)} {\mu(B^o(x,\Delta/4)\cap Z)} \text{\quad and\quad } g(u)=\mu(B(u,\diam(\bbQ)/4)\cap \bbQ) 
\]
reach  maxima on the  compact domains $Z$ and $\bbQ$ (respectively).
While $f$ and $g$ are not necessarily continuous, the assumption is indeed correct. 
To see this for $f$ (the argument for $g$ is similar),
observe first that $0\le f(x)\le 1$.  
Let $(a_i)_i$ be a sequence of points in $Z$ 
such that $f(a_i)\nearrow \sup_{x\in Z} f(x)$.
Since $Z$ is compact, by moving to a subsequence we may assume that $(a_i)_i$  converges to some point $\lim_i a_i= a\in Z$, and by moving further to a subsequence, that $d(a_i,a)\searrow0$.
Next we claim that 
$\varlimsup_i \mu(B_d(a_i, \Delta/8))\le \mu(B_d(a,\Delta/8))$.
This is proved by observing that
\(
B_d(a_i, \Delta/8) \subseteq B_d(a,\Delta/8+d(a_i,a)),
\)
and therefore
\begin{multline*}
\mu(B_d(a_i, \Delta/8)) - \mu(B_d(a,\Delta/8))
\\ \le \mu(B_d(a_, \Delta/8+d(a_i,a))) - \mu(B_d(a,\Delta/8))
= \mu \bigl(B_d(a, \Delta/8+d(a_i,a))\setminus 
B_d(a,\Delta/8)\bigr).
\end{multline*}
The right-hand side converges to $0$, since the sequence of subsets decreases to $\emptyset$.

In a similar fashion,
\begin{multline*}
\mu(B^o_d(a_i, \Delta/4)) - \mu(B^o_d(a,\Delta/4))
\\ \ge \mu(B^o_d(a, \Delta/4-d(a_i,a))) - \mu(B^o_d(a,\Delta/4))
= -\mu \bigl(B^o_d(a, \Delta/4)\setminus 
B^o_d(a,\Delta/4-d(a_i,a))\bigr),
\end{multline*}
and again the right-hand side converges to $0$, since the sequence of subsets decreases to $\emptyset$.
Hence,
$\varliminf_i \mu(B^o_d(a_i, \Delta/4))\ge \mu(B^o_d(a,\Delta/4 ))$. We summarize:
\[f(a) \le \sup_{x\in Z}f(x)=
\lim_i f(a_i)\le 
\frac{\varlimsup_i \mu(B_d(a_i, \Delta/8))}{\varliminf_i \mu(B^o_d(a_i, \Delta/4))}
\le \frac{\mu(B_d(a,\Delta/8)}{\mu(B^o_d(a,\Delta/4 )}
= f(a). 
\]
\end{remark}

\begin{proof}[Proof of Lemma~\ref{lem:compact-um}]
That $\rho$ is an ultrametric is a straight-forward conclusion
from the order on trees: $x\wedge z \succeq (x\wedge y)\wedge (y\wedge z)$, and therefore 
\[
\rho(x,z)=\Delta(x\wedge z)\le \Delta (x\wedge y)\wedge (y\wedge z)) \le \max\{\Delta(x\wedge y), \Delta (y\wedge z)\}
=\max\{\rho(x,y),\rho(y,z) \}.
\]

Next we prove that $\mathcal O_\T$ is the set of open balls.
Let $v\in\T$ such that $v=\emptyset$ or $\Delta(v)<\Delta(u)$, where $u$ is the parent of $v$. Let
$x\in\partial_\T(v)$ be an arbitrary branch in the boundary of $v$, and $r\in(\Delta(v),\Delta(u))$. Then clearly,
$\partial_\T(v)=B_\rho^o(x,r)$.
In the other direction, fix an open ball $B_\rho^o(x,r)$, where $x\in\partial_\T(\emptyset)$, and $r>0$. 
Since $x\in\partial_\T(\emptyset)$ is a branch of $\T$, let $v\in x$ a vertex on the branch that satisfies $\Delta(v)<r$, and for any $u\prec v$, $r \le \Delta(u)$.
Clearly, such $v$ exists and is unique.
Then for any $y\in \partial_\T(v)$, $x\wedge y \succeq v$, and hence $\rho(x,y)<r$, and on the other hand for any $y$ for which $
\rho(x,y)<r$, $x\wedge y\in x$ and $x \wedge y \succ v$, hence
$y\in \partial_\T(v)$. This proves that $B^o(x,r)=\partial_\T(v)$.
Similar arguments, but replacing strict inequalities of $r$ with weak inequalities, prove that $\mathcal O_\T$ is also the set of closed balls with positive radii in $(\partial_\T(\emptyset)),\rho)$ .

We are left to prove compactness of $\partial_\T(\emptyset)$. 
Assume $\partial_\T(\emptyset)$ is infinite (otherwise, it is trivially compact). 
Fix an infinite $A\subseteq \partial_\T(\emptyset)$. We should prove that $A$ has an accumulation point in $\partial_\T(\emptyset)$. 
To achieve it, we construct an infinite non-increasing sequence of infinite subsets $A=A_0\supseteq A_1\supseteq A_2\supseteq \ldots$ and an infinite  sequence of vertices $\emptyset=v_0,v_1, v_2,\ldots$ in $\T$ satisfying 
$A_i\subseteq \partial_\T(v_i)$, and 
$v_{i+1}\in \{v_{i}{}^\smallfrown 0, v_{i}{}^\smallfrown 1\}$. 
The construction is by induction.
In the base case $v_0=\emptyset$, $A=A_0\subseteq\partial_\T(\emptyset)$ holds by assumption. 
Assume we have already defined $A_0,\ldots, A_n$ and $v_0,\ldots,v_n$ to satisfy the above when $i\in\{0,\ldots,n-1\}$. 
since $A_n\subset \partial_\T(v_n)$ is infinite, $v_n$ is not a leaf and hence have two children, and
\[
A_n=A_n\cap \partial(v_n)=
(A_n\cap \partial_\T(v_n{}^\smallfrown 0)) \cup 
(A_n\cap \partial_\T(v_n{}^\smallfrown 1)).
\]
Hence, at least one of 
$A_n\cap \partial_\T(v_n{}^\smallfrown 0)$ and
$A_n\cap \partial_\T(v_n{}^\smallfrown 1)$ 
must be infinite.
Let $v_{n+1}\in \{v_n{}^\smallfrown 0, v_n{}^\smallfrown 1\}$
for which 
$A_n\cap \partial_\T(v_{n+1})$ is infinite, and define
$A_{n+1}=A_n\cap \partial_\T(v_{n+1})$.
The sequence $b=(v_n)_n$ is 
a  branch in $\partial_\T(u)$.
For every $n\in\mathbb N$ we have, by the covering property,
\[
\emptyset\ne A_n\subseteq A\cap \partial_\T(v_n)\subseteq B_\rho(b,\Delta(v_n)).
\]
Hence, $B_\rho({b},\Delta(v_n))\cap A \supseteq A_n\ne \emptyset$.  Since $\Delta(v_n)\to 0$, it means that 
${b}\in\partial_\T(u)$ is an accu\-mu\-lation point of~$A$.
\end{proof}

\subsection*{Acknowledgments}
The author thanks an anonymous referee for insightful comments and ideas. 
In particular, Section~\ref{sec:remarks} and Section~\ref{sec:zindulka}
are based on --- and are an outgrowth of --- suggestions of the referee.
The author was supported by BSF grant no{.} 2018223.

The author states no conflict of  interest.

\bibliographystyle{plainurl}
\bibliography{simple-ramsey} 

\end{document}